\documentclass[10pt]{article}
\usepackage{cite}
\usepackage{mathrsfs}
\usepackage{amsfonts}
\usepackage{amsmath}
\usepackage{amsfonts,amssymb}
\usepackage{dsfont}
\usepackage{curves}
\usepackage{mathrsfs}
\usepackage{pifont}
\usepackage{amssymb}
\allowdisplaybreaks

\numberwithin{equation}{section}

\date{}

\def\BigRoman{\uppercase\expandafter{\romannumeral\number\count 255 }}
\def\Romannumeral{\afterassignment\BigRoman\count255=}

\begin{document}
\title{Toughness and spectral radius in graphs
}
\author{\small Sufang Wang$^{1}$\footnote{E-mail address: wangsufangjust@163.com (S. Wang).},
Wei Zhang$^{2}$\footnote{Corresponding author. E-mail address: zw\_wzu@163.com (W. Zhang).}\\
\small $1$.  School of Public Management, Jiangsu University of Science and Technology,\\
\small Zhenjiang, Jiangsu 212100, China\\
\small $2$.  School of Economics and Management, Wenzhou University of Technology,\\
\small Wenzhou, Zhejiang 325000, China\\
}

\maketitle
\begin{abstract}
\noindent Let $t$ be a positive integer, and let $G$ be a connected graph of order $n$ with $n\geq t+2$. A graph $G$ is said to be
$\frac{1}{t}$-tough if $|S|\geq\frac{1}{t}c(G-S)$ for every subset $S$ of $V(G)$ with $c(G-S)\geq2$, where $c(G-S)$ is the number of connected
components in $G-S$. The adjacency matrix of $G$ is denoted by $A(G)$. Let $\lambda_1(G)\geq\lambda_2(G)\geq\dots\geq\lambda_n(G)$ be the eigenvalues
of $A(G)$. In particular, the eigenvalue $\lambda_1(G)$ is called the spectral radius of $G$. In this paper, we prove that $G$ is a
$\frac{1}{t}$-tough graph unless $G=K_1\vee(K_{n-t-1}\cup tK_1)$ if $\lambda_1(G)\geq\eta(t,n)$, where $\eta(t,n)$ is the largest root of
$x^{3}-(n-t-2)x^{2}-(n-1)x+t(n-t-2)=0$.
\\
\begin{flushleft}

(2020) Mathematics Subject Classification: 05C50, 90B99
\end{flushleft}
\end{abstract}

\section{Introduction}

In this paper, we deal with simple and undirected graphs. Let $G$ be a graph with vertex set $V(G)$ and edge set $E(G)$. The order of $G$
is denoted by $|V(G)|=n$. For any $S\subseteq V(G)$, we denote by $G[S]$ the subgraph of $G$ induced by $S$, and by $G-S$ the subgraph
formed from $G$ by deleting the vertices in $S$ and their incident edges. Let $K_n$ denote a complete graph of order $n$. For two
vertex-disjoint graphs $G_1$ and $G_2$, we use $G_1\cup G_2$ to denote the disjoint union of $G_1$ and $G_2$. The join $G_1\vee G_2$ is
the graph obtained from $G_1\cup G_2$ by adding all possible edges between $V(G_1)$ and $V(G_2)$.

Let $V(G)=\{v_1,v_2,\ldots,v_n\}$. The adjacency matrix of $G$ is denoted by $A(G)$. Let $\lambda_1(G)\geq\lambda_2(G)\geq\dots\geq\lambda_n(G)$
be the eigenvalues of $A(G)$. In particular, the eigenvalue $\lambda_1(G)$ is called the spectral radius of $G$.

Let $t$ be a positive real number. A graph $G$ is $t$-tough if $tc(G-S)\leq|S|$ for every subset $S\subseteq V(G)$ with $c(G-S)\geq2$,
where $c(G-S)$ is the number of components in $G-S$. The toughness $t(G)$ of $G$ is the largest value of $t$ for which $G$ is $t$-tough
(taking $t(K_n)=\infty$). Hence, if $G$ is not a complete graph, $t(G)=\min\{\frac{|S|}{c(G-S)}: S\subseteq V(G), c(G-S)\geq2\}$, where
the minimum is taken over all cut set of vertices in $G$.

Many authors \cite{N,FN,LLD,O,OPPZ,ZLt,ZLc,ZSL1,ZSL2,ZZS,ZZ,ZZL} investigated some interesting spectral properties of $A(G)$. Some properties on toughness of
graph were investigated in \cite{GJS,K,YH,KKV,Zhs,ZWX,ZWB,ZSB,ZSY,ZPX,WZi}. The relationship between toughness and eigenvalues was first studied by
Alon \cite{A}. Brouwer \cite{B} independently established a relationship between toughness and eigenvalues. More results on the relationship
between toughness and eigenvalues can be found in \cite{Gt,Ga,CW}. Fan, Lin and Lu \cite{FLLt} provided two spectral radius conditions
for a graph to be $t$-tough. In this paper, we continue to investigate the existence of a tough graph, and obtain a new sufficient
condition for a graph to be $\frac{1}{t}$-tough by using spectral radius.

\medskip

\noindent{\textbf{Theorem 1.1.}} Let $t$ be a positive integer, and let $G$ be a connected graph of order $n\geq t+2$. If
$$
\lambda_1(G)\geq\eta(t,n),
$$
then $G$ is a $\frac{1}{t}$-tough graph unless $G=K_1\vee(K_{n-t-1}\cup tK_1)$, where $\eta(t,n)$ is the largest root of
$x^{3}-(n-t-2)x^{2}-(n-1)x+t(n-t-2)=0$.

\medskip

We have the following corollary if $t=1$ in Theorem 1.1.

\medskip

\noindent{\textbf{Corollary 1.2.}} Let $G$ be a connected graph of order $n\geq3$. If
$$
\lambda_1(G)\geq\eta(n),
$$
then $G$ is a 1-tough graph unless $G=K_1\vee(K_{n-2}\cup K_1)$, where $\eta(n)$ is the largest root of $x^{3}-(n-3)x^{2}-(n-1)x+n-3=0$.

\section{Preliminaries}

In this section, we put forward some necessary preliminary lemmas, which will be used to prove our main results.

\medskip

\noindent{\textbf{Lemma 2.1}} (Li and Feng \cite{LF}). If $G$ is a connected graph, and $H$ is a proper subgraph of $G$, then
$$
\lambda_1(G)>\lambda_1(H).
$$

\medskip

Let $M$ be a real matrix whose rows and columns are indexed by $V=\{1,2,\ldots,n\}$. Assume that $M$, with respect to the partition
$\pi:V=V_1\cup V_2\cup\cdots\cup V_m$, can be written as
\begin{align*}
M=\left(
  \begin{array}{cccc}
    M_{11} & M_{12} & \cdots & M_{1m}\\
    M_{21} & M_{22} & \cdots & M_{2m}\\
    \vdots & \vdots & \ddots & \vdots\\
    M_{m1} & M_{m2} & \cdots & M_{mm}\\
  \end{array}
\right),
\end{align*}
where $M_{ij}$ denotes the submatrix (block) of $M$ formed by rows in $V_i$ and columns in $V_j$. Let $q_{ij}$ denote the average row sum of
$M_{ij}$, namely, $q_{ij}$ is the sum of all entries in $M_{ij}$ divided by the number of rows. Then matrix $M_{\pi}=(q_{ij})$ is called the
quotient matrix of $M$. If the row sum of every block $M_{ij}$ is a constant, then the partition is equitable.

\medskip

\noindent{\textbf{Lemma 2.2}} (You, Yang, So and Xi \cite{YYSX}). Let $M$ be a real matrix with an equitable partition $\pi$,
and let $M_{\pi}$ be the corresponding quotient matrix. Then each eigenvalue of $M_{\pi}$ is an eigenvalue of $M$. Furthermore, if $M$ is
nonnegative, then the spectral radiuses of $M$ and $M_{\pi}$ are equal.

\medskip

\noindent{\textbf{Lemma 2.3}} (Fan, Goryainov, Huang and Lin \cite{FGHL}). Let $n_1\geq n_2\geq\cdots\geq n_c$ be positive integers with
$n=\sum_{i=1}^{c}n_i+s$ and $n_1\leq n-s-c+1$. Then
$$
\lambda_1(K_s\vee(K_{n_1}\cup K_{n_2}\cup\cdots\cup K_{n_c}))\leq\lambda_1(K_s\vee(K_{n-s-c+1}\cup(c-1)K_1)),
$$
where the equality holds if and only if $(n_1,n_2,\ldots,n_c)=(n-s-c+1,1,\ldots,1)$.

\section{The proof of Theorem 1.1}

In this section, we are to prove Theorem 1.1, which gives a spectral radius condition for the existence of a $\frac{1}{t}$-tough graph.

\medskip

\noindent{\it Proof of Theorem 1.1.} Let $\Phi(x)=x^{3}-(n-t-2)x^{2}-(n-1)x+t(n-t-2)$ be a real function in $x$ and let $\eta(t,n)$ be the
largest root of $\Phi(x)=0$. Suppose that $G$ is not a $\frac{1}{t}$-tough graph, there exists some nonempty subset $S$ of $V(G)$ satisfying
$c(G-S)\geq t|S|+1$. Let $|S|=s$. Then $G$ is a spanning subgraph of $G_1=K_s\vee(K_{n_1}\cup K_{n_2}\cup\cdots\cup K_{n_{ts+1}})$ for
some positive integers $n_1\geq n_2\geq\cdots\geq n_{ts+1}$ with $\sum\limits_{i=1}^{ts+1}n_i=n-s$. Applying Lemma 2.1, we get
\begin{align}\label{eq:3.1}
\lambda_1(G)\leq\lambda_1(G_1),
\end{align}
where the equality holds if and only if $G=G_1$.

Let $G_2=K_s\vee(K_{n-ts-s}\cup tsK_1)$. In terms of Lemma 2.3, we obtain
\begin{align}\label{eq:3.2}
\lambda_1(G_1)\leq\lambda_1(G_2),
\end{align}
where the equality holds if and only if $(n_1,n_2,\ldots,n_{ts+1})=(n-ts-s,1,\ldots,1)$. We are to prove the following claim.

\medskip

\noindent{\bf Claim 1.} $\lambda_1(G_2)\leq\eta(t,n)$ with equality if and only if $G_2=K_1\vee(K_{n-t-1}\cup tK_1)$.

\noindent{\it Proof.} By virtue of the partition $V(G_2)=V(K_s)\cup V(K_{n-ts-s})\cup V(tsK_1)$, the quotient matrix of $A(G_2)$ is
\begin{align*}
B_1=\left(
  \begin{array}{ccc}
    s-1 & n-ts-s & ts\\
    s & n-ts-s-1 & 0\\
    s & 0 & 0\\
  \end{array}
\right).
\end{align*}
Then the characteristic polynomial of $B_1$ is
$$
\Phi_{B_1}(x)=x^{3}-(n-ts-2)x^{2}-(n+ts^{2}-ts-1)x+ts^{2}(n-ts-s-1).
$$
Since the partition $V(G_2)=V(K_s)\cup V(K_{n-ts-s})\cup V(tsK_1)$ is equitable, in view of Lemma 2.2, the largest root, say $\eta_1$, of
$\Phi_{B_1}(x)=0$ is equal to $\lambda_1(G_2)$.

Note that $n-ts-s\geq1$ and $K_s\vee(ts+1)K_1$ is a subgraph of $G_2$. The quotient matrix of $A(K_s\vee(ts+1)K_1)$ with respect to the partition
$V(K_s\vee(ts+1)K_1)=V(K_s)\cup V((ts+1)K_1)$ is
\begin{align*}
B_2=\left(
  \begin{array}{ccc}
    s-1 & ts+1\\
    s & 0\\
  \end{array}
\right).
\end{align*}
Then the characteristic polynomial of $B_2$ is
$$
\Phi_{B_2}(x)=x^{2}-(s-1)x-s(ts+1).
$$

Note that the partition $V(K_s\vee(ts+1)K_1)=V(K_s)\cup V((ts+1)K_1)$ is equitable. By virtue of Lemma 2.2, the largest root, say $\eta_2$,
of $\Phi_{B_2}(x)=0$ equals $\lambda_1(K_s\vee(ts+1)K_1)$. According to the root formula, we conclude
\begin{align}\label{eq:3.3}
\lambda_1(K_s\vee(ts+1)K_1)=\eta_2=\frac{s-1+\sqrt{(4t+1)s^{2}+2s+1}}{2}.
\end{align}
According to \eqref{eq:3.3} and Lemma 2.1, we obtain
\begin{align}\label{eq:3.4}
\eta_1=\lambda_1(G_2)\geq\lambda_1(K_s\vee(ts+1)K_1)=\eta_2=\frac{s-1+\sqrt{(4t+1)s^{2}+2s+1}}{2}.
\end{align}

Note that $\Phi_{B_1}(\eta_1)=0$. By a simple calculation, we get
\begin{align}\label{eq:3.5}
\Phi(\eta_1)=\Phi(\eta_1)-\Phi_{B_1}(\eta_1)=(s-1)h_1(\eta_1),
\end{align}
where $h_1(\eta_1)=-t\eta_1^{2}+ts\eta_1-(ts+t)n+t^{2}s^{2}+ts^{2}+t^{2}s+2ts+t^{2}+2t$. If $s=1$, then by \eqref{eq:3.5} we have $\Phi(\eta_1)=0$.
Thus $\eta_1$ is a root of $\Phi(x)=0$ and so
\begin{align}\label{eq:3.6}
\lambda_1(G_2)=\eta_1\leq\eta(t,n),
\end{align}
where the last equality holds if and only if $G_2=K_1\vee(K_{n-t-1}\cup tK_1)$. Next, we consider $s\geq2$.

Recall that $n\geq ts+s+1$. Then
\begin{align}\label{eq:3.7}
h_1(\eta_1)=&-t\eta_1^{2}+ts\eta_1-(ts+t)n+t^{2}s^{2}+ts^{2}+t^{2}s+2ts+t^{2}+2t\nonumber\\
\leq&-t\eta_1^{2}+ts\eta_1-(ts+t)(ts+s+1)+t^{2}s^{2}+ts^{2}+t^{2}s+2ts+t^{2}+2t\nonumber\\
=&-t\eta_1^{2}+ts\eta_1+t^{2}+t.
\end{align}
Let $f_1(x)=-tx^{2}+tsx+t^{2}+t$ be a real function in $x$. Then the symmetry axis of $f_1(x)$ is $x=\frac{s}{2}$. Clearly, $f_1(x)$ is
decreasing when $x\geq\frac{s}{2}$. Combining this with \eqref{eq:3.4}, we obtain
\begin{align*}
f_1(\eta_1)\leq&f_1\left(\frac{s-1+\sqrt{(4t+1)s^{2}+2s+1}}{2}\right)\\
=&\frac{t}{2}(-2ts^{2}+\sqrt{(4t+1)s^{2}+2s+1}-s+2t+1).
\end{align*}
Note that
$$
(4t+1)s^{2}+2s+1-(t+2)^{2}s^{2}=-(t^{2}+3)s^{2}+2s+1<0.
$$
Thus, we conclude
\begin{align*}
f_1(\eta_1)\leq&\frac{t}{2}(-2ts^{2}+\sqrt{(4t+1)s^{2}+2s+1}-s+2t+1)\\
<&\frac{t}{2}(-2ts^{2}+(t+2)s-s+2t+1)\\
=&\frac{t}{2}(-2ts^{2}+(t+1)s+2t+1).
\end{align*}
Note that
$$
\frac{t+1}{4t}<2\leq s.
$$
Then
\begin{align*}
f_1(\eta_1)<&\frac{t}{2}(-2ts^{2}+(t+1)s+2t+1)\\
\leq&\frac{t}{2}(-8t+2(t+1)+2t+1)\\
=&\frac{t}{2}(-4t+3)\\
<&0.
\end{align*}
Combining this with \eqref{eq:3.5}, \eqref{eq:3.7} and $s\geq2$, we obtain
$$
\Phi(\eta_1)=(s-1)h_1(\eta_1)\leq(s-1)f_1(\eta_1)<0,
$$
which yields that
\begin{align}\label{eq:3.8}
\lambda_1(G_2)=\eta_1<\eta(t,n).
\end{align}

In terms of \eqref{eq:3.6} and \eqref{eq:3.8}, we conclude
$$
\lambda_1(G_2)\leq\eta(t,n),
$$
where the equality holds if and only if $G_2=K_1\vee(K_{n-t-1}\cup tK_1)$. Claim 1 is proved. \hfill $\Box$

It follows from \eqref{eq:3.1}, \eqref{eq:3.2} and Claim 1 that
$$
\lambda_1(G)\leq\eta(t,n),
$$
with equality if and only if $G=K_1\vee(K_{n-t-1}\cup tK_1)$. Which is a contradiction to the spectral radius condition of Theorem 1.1. Also note
that $\lambda_1(K_1\vee(K_{n-t-1}\cup tK_1))$ is the largest root of $x^{3}-(n-t-2)x^{2}-(n-1)x+t(n-t-2)=0$. This completes the proof of
Theorem 1.1. \hfill $\Box$

\section*{Data availability statement}

My manuscript has no associated data.

\section*{Declaration of competing interest}

The authors declare that they have no conflicts of interest to this work.

\section*{Contributions}

Sufang Wang Writing-supervision, conceptualization, methodology, original draft. Wei Zhang Methodology, formal analysis.


\end{document}